\documentclass[11pt,twoside]{article}

\usepackage{fancyhdr}
\usepackage{epsfig}
\usepackage{latexsym}
\usepackage{microtype}
\usepackage{lastpage} 
\usepackage{hyperref}
\usepackage{graphicx}
\usepackage{lipsum}
\usepackage{amsmath}
\usepackage{lineno,hyperref}
\usepackage{mathtools}
\usepackage{amssymb}
\usepackage{multirow}
\usepackage{comment}
\usepackage[dvipsnames*,svgnames]{xcolor}
\hypersetup{colorlinks=true,allcolors=DarkBlue,linkcolor=black}

\DeclarePairedDelimiter\ceil{\lceil}{\rceil}
\DeclarePairedDelimiter\floor{\lfloor}{\rfloor}

\newtheorem{corollary}{Corollary}
\newtheorem{lemma}{Lemma}

\newcommand{\proof}[1]{\par\noindent {\bf Proof:}\ \ \ #1 \hfill$\Box$\hspace{1ex}}

\newcommand{\trshortyear}{21}
\newcommand{\trpapernumber}{01}
\newcommand{\trmonth}{June}
\newcommand{\tryear}{2021}

\newcommand{\TheAuthor}{}
\newcommand{\Author}[1]{\renewcommand{\TheAuthor}{#1}}

\newcommand{\TheTitle}{}
\newcommand{\Title}[1]{\renewcommand{\TheTitle}{#1}}
\Author{E. F. OLARIU, C. FR\u ASINARU}
\Title{Improving Lower Bounds for Equitable Chromatic Number}

\newcommand\blfootnote[1]{%
	\begingroup
	\renewcommand\thefootnote{}\footnote{#1}%
	\addtocounter{footnote}{-1}%
	\endgroup
}

\begin{document}

	\blfootnote{The publisher does not claim any copyright for the technical reports. The author keeps the full copyright for the paper, and is thus free to transfer the copyright to a publisher if the paper is accepted for publication elsewhere. }   	

\parindent=8mm

\noindent {{\bf \scriptsize Faculty of Computer Science, Alexandru Ioan Cuza University Ia\c si}}

\noindent {{\bf \scriptsize Technical Report TR \trshortyear-\trpapernumber, \trmonth ~ \tryear}}
\vskip -3mm
\noindent\rule{10.2cm}{0.4pt}
\vskip -1mm
\noindent

\vspace{1cm}
\begin{center}
{\Large\bf Improving lower bounds \\
for equitable chromatic number
}
\end{center}
\vspace{4mm}

\begin{center}
{\large Emanuel Florentin OLARIU, Cristian FR\u ASINARU}\\
Faculty of Computer Science, Alexandru Ioan Cuza University Ia\c si, \\
General Berthelot 16, 700483 Ias\c i, Romania, \\
Email: {\tt olariu@info.uaic.ro, acf@info.uaic.ro}
\end{center}
\vspace{3ex}

\date{}

\begin{abstract}     
	
In many practical applications the underlying graph must be as equitable colored as possible. A coloring is called {\it equitable} if the number of vertices colored with each color differs by at most one, and the least number of colors for which a graph has such a coloring is called its equitable chromatic number.

We introduce a new integer linear programming approach for studying the equitable coloring number of a graph and show how to use it for improving lower bounds for this number. The two stage method is based on finding or upper bounding the maximum cardinality of an equitable color class in a valid equitable coloring and, then, sequentially improving the lower bound for the equitable coloring number.

The computational experiments were carried out on DIMACS graphs and other graphs from the literature.

\smallskip


\end{abstract}

\section{Introduction}                                 

In graph theory there exists a wide range of optimization problems with pertinent practical importance, and one of the most studied is the Graph Coloring Problem (GCP); this problem arises in many applications such as scheduling, timetabling, electronic bandwidth allocation and sequencing problems (see \cite{malaguti10} for a survey).

Given $ G = (V, E)  $ a graph, where $ V $ is the set of vertices and $ E $ is the set of edges, a {\it $ p $-coloring} of (the vertices of) $ G $ is a map $ c : V \to \{ 1, 2, \ldots, p \} $ such that any two adjacent vertices have different colors. Vertices with the same color make together a color class: $ c^{-1}(i) $, for $ 1 \le i \le p $; some of these color classes could be empty but all of them are stable sets of $ G $. The {\it graph coloring problem} for $ G $ consists of finding the minimum number $ p $ such that~$ G $ has a $ p $-coloring. This minimum number of colors is called the {\it chromatic number} of the graph G and is denoted by $ \chi(G) $.

One of the usual applications of this problem is to model the simple scheduling problem: assign a given set of tasks to workers, knowing that some pairs of tasks cannot be assigned to the same worker. We can model this assignment problem by building a graph whose vertices are the tasks and whose edges are the conflicting pairs of tasks. A coloring of the resulted graph will be a feasible assignment of all tasks while the chromatic number will be the minimum number of needed workers.

By imposing additional restrictions one may get variations of the graph coloring problem. For instance, in the above scheduling problem, it may be required to ensure a certain kind of load balancing of workload. This can be viewed as an {\it equity constraint} and by imposing it we get the {\it equitable coloring problem}.

An {\it equitable $ p $-coloring} of $ G $ is $ p $-coloring such that difference on the cardinalities of any two non-empty color classes is at most one. Each subset is associated with a color and called a color set. The Equitable Coloring Problem (ECP) consists of finding the minimum value $ p $ such that there is an equitable $ p $-coloring of G. This value is said to be the equitable chromatic number of G and is denoted by $ \chi_=(G) $ or $ \chi_{eq}(G) $.

The equitable coloring problem was introduced by Meyer in \cite{meyer73} motivated by some scheduling problems; other applications of this problem are partitioning and load balancing in multiprocessor machines \cite{blazewicz01}, scheduling~\cite{irani96}, probability theory, municipal garbage collection \cite{tucker73}. A review of some basic results on ECP are provided in \cite{furmanczyk04,furmanczyk16}.

Computing the equitable chromatic number was proved to be NP-hard (see \cite{furmanczyk05}) and the number of graph families for which ECP is  known to be easy to solve is small (the trees, the complete $ n $-partite graphs, the wheel graphs, the graphs with bounded treewidth etc). Several exact and heuristic
approaches are known in the literature for solving the ECP for arbitrary graphs: tabu-search heuristics (\cite{lai15,mendez14a}), linear programming based algorithms (\cite{bahiense14,mendez14b}), degree of saturation based heuristic (\cite{mendez15}), greedy-based constructive heuristics (\cite{furmanczyk04}) etc.

In this paper we introduce a new method for finding lower bounds for the equitable chromatic number. Our approach starts by finding an upper bound of the maximum cardinality of an equitable color class and, then, improves the lower bound on the equitable chromatic number by verifying one by one the consecutive possible values.

In order to do this we need two integer linear programming models: one for finding upper bounds of the maximum cardinality of an equitable color class and one for deciding if for a given $ p $, the graph admits an equitable $ p $-coloring. Both these models are based on integer linear programming models for the classic coloring problem using partial-ordering and assignment models. 

The remaining of the paper is organized as follows: in Section \ref{section:background} we discuss our setting and the related work, in Section \ref{section:LPpartial} we describe and explain the LP models, Section \ref{section:numerical} contains the numerical results, and the last section is dedicated to conclusions.

%Insert here the text of your paper. \\
%Cite references in the normal manner: \cite{Babuska}, \cite{Brauner2010},  and 
%\cite{ST99}.  \\
%Please provide {\tt doi} information in your BibTeX entries.

\section{Background}\label{section:background}

Let $ G = (V, E) $ be an undirected graph with $ n $ vertices. While a graph admitting a $ p $-coloring admits also a $ (p + 1) $-coloring, this is not necessarily true for equitable colorings. Other differences between the two numbers, $ \chi $ and $ \chi_{eq} $: the equitable chromatic number of a subgraph is not necessarily smaller than the equitable chromatic number of the main graph, and the equitable chromatic numbers of connected components are not related with the  equitable chromatic number of the graph itself.

A common property of the two numbers is the following $ \chi(G), \chi_{eq}(G) $ are both at most $ \Delta(G) + 1 $ - the second inequality was obtained by Hajnal and Szemeredi:  (\cite{hajnal70,kierstead10}), while the first is due to the greedy coloring algorithm. A simple observation shows that if $ G $ admits an equitable $ p $-coloring whose non-empty color classes are $ S_1, S_2, \ldots S_p $, then $ |S_i| \in \{ \displaystyle \lfloor n/p \rfloor, \lceil n/p \rceil \} $, $ i = \overline{1, p} $. 

There are three main linear programming approaches to graph coloring problem resolution: the {\it assignment} model, the {\it representatives} model, and the {\it set covering} model. The straightforward way of modeling the equitable coloring problem is the assignment model that uses two types of variables:~$ x_{vi} $, with $ v \in V $ and $ i \in \{1, 2, \ldots, n \} $, $ x_{vi} = 1 $ if and only if vertex $ v $ receives the color $ i $, and $ w_i $ with $ i \in \{1, 2, \ldots, n \} $, where $ w_i = 1 $ if and only if color $ i $ is~used.
\begin{align}
min & \begin{array}{l}
  \displaystyle  \sum_{i = 1}^n w_i 
\end{array}
\notag\\
& \begin{array}{l}
\displaystyle \sum_{i = 1}^n x_{vi}=  1, \forall v \in V,
\end{array}
\label{asgn1}\\
& \begin{array}{l}
\displaystyle x_{ui} + x_{vi} \le  w_i, \forall uv \in E, i = \overline{1, n}
\end{array}
\label{asgn2}\\
& \begin{array}{l}
\displaystyle x_{vi} \in \{ 0, 1 \} \;  \forall v \in V, i = \overline{1, k}, \; w_i \in \{ 0, 1 \} \; i = \overline{1, n}
\end{array}
\notag
\end{align}

\noindent To this classical model some other constraints can be added in order to remove symmetric solutions (see \cite{mendez08}); for obtaining $ p $-equitable colorings we add the following constraints (see \cite{mendez14b}):
\begin{align}
& \begin{array}{l}
\displaystyle \sum_{i = 1}^n x_{vi} \ge \sum_{j = i}^n  \displaystyle \floor*{\frac{n}{p}} (w_j - w_{j + 1}), \; i = \overline{1, n - 1},
\end{array}
\label{asgn3}\\
& \begin{array}{l}
\displaystyle \sum_{i = 1}^n x_{vi} \le \sum_{j = i}^n  \displaystyle \ceil*{\frac{n}{p}} (w_j - w_{j + 1}), \; i = \overline{1, n - 1},
\end{array}
\label{asgn4}\\
& \begin{array}{l}
\displaystyle x_{ui} \le  w_i, \forall u \mbox{ isolated}, i = \overline{1, n}
\end{array}
\label{asgn5}
\end{align}

\noindent In the representatives model  (see \cite{jabrayilov18}) each color class is represented by exactly one vertex. The binary variables are: $ x_{uv} $, with $ u, v \in V $ such that $ uv \notin E $ (including here the case $ u = v $); for $ u \not= v$, $ x_{uv} = 1 $ if and only if $ v $ is represented by $ u $, while for $ u = v $, $ x_{vv} = 1 $ if and only if $ v $ is representative for its color class.
\begin{align}
min & \begin{array}{l}
  \displaystyle  \sum_{v \in V} x_{vv}
\end{array}
\notag\\
& \begin{array}{l}
\displaystyle \sum_{u: uv \notin E}x_{uv} \ge  1, \forall v \in V,
\end{array}
\label{rep1}\\
& \begin{array}{l}
\displaystyle x_{uv} + x_{vw} \le  x_{vv}, \forall v \in V, \forall u, w \in V\setminus \{ v \}  \mbox{ s. t. } uv, vw \notin E,
\end{array}
\label{rep2}\\
& \begin{array}{l}
\displaystyle x_{uv} \in \{ 0, 1 \} \;  \forall uv \notin E
\end{array}
\notag
\end{align}

\noindent Some other constraints must be added for obtaining an equitable coloring of minimum number of colors. The difference from the above model is that by solving this modified model (see \cite{bahiense14}) we get the equitable chromatic number  -- $ \chi_{eq}(G) $.

The set covering based model \cite{mehrotra96} aims to cover all vertices with the minimum number of stable (vertex independent) sets.
\begin{align}
min & \begin{array}{l}
  \displaystyle  \sum_{S \in {\mathcal S}} x_S
\end{array}
\notag\\
& \begin{array}{l}
\displaystyle \sum_{S \in {\mathcal S}: v \in S} x_S \ge  1, \forall v \in V,
\end{array}
\label{set1}\\
& \begin{array}{l}
\displaystyle x_S \in \{ 0, 1 \} \;  \forall S \in {\mathcal S}
\end{array}
\notag
\end{align}
$ {\mathcal S} $ is the family of all stable sets for the graph coloring problem, but in order to get a $ p $-equitable coloring we must change this family to $ {\mathcal S'} $:
\begin{align}
& \begin{array}{l}
\displaystyle |S| \le \displaystyle \ceil*{\frac{n}{p}}, \forall S \in {\mathcal S'},
\end{array}
\label{set2}\\
& \begin{array}{l}
\displaystyle |S| \ge  \displaystyle \floor*{\frac{n}{p}}, \forall S \in {\mathcal S'},
\end{array}
\label{set3}
\end{align}

\noindent For finding a $ p $-equitable coloring (if any) the above model must be solved by the column generation method combined with other approaches (e. g. branch-and-bound giving a branch-and-price algorithm).

For our approach we will investigate only the assignment and the partial-ordering based models.

\section{LP Partial-Ordering and\\ Assignment Based Models}\label{section:LPpartial}

\subsection{LP Partial-Ordering Model}
In the partial-ordering based integer linear programming model for graph coloring problem (see \cite{jabrayilov18}) the colors are linearly ordered and they are not directly assigned to vertices, but we determine a relative order of each vertex with respect to each color. If $ v $ is a vertex and $ i $ is a color, a solution to this model establishes that if $ v $ is neither greater nor lower than $ i $, then $ v $ will receive the color $ i $.  

The (binary) variables of the model are: $ y_{i, v} $ and $ z_{v, i} $ with $ i \in \{1, 2, \ldots, k \} $ and $ v \in V $, where $ k $ is an upper bound for the chromatic number; $ y_{i,v} = 1 $ if and only if $ v $ is greater than $ i $, while $ z_{v, i} = 1 $ if and only if $ v $ is less than $ i $
\begin{align}
min & \begin{array}{l}
 \displaystyle 1+ \left( \sum_{i = 1}^k y_{i, q} \right)
\end{array}
\notag\\
& \begin{array}{l}
\displaystyle z_{v, 1} =  0, \forall v \in V,
\end{array}
\label{order1}\\
& \begin{array}{l}
\displaystyle y_{k, v} = 0, \forall v \in V,
\end{array}
\label{order2}\\
& \begin{array}{l}
\displaystyle y_{i, v} - y_{i + 1, v} \ge  0, \forall v \in V, i = \overline{1, k - 1},
\end{array}
\label{order3}\\
& \begin{array}{l}
\displaystyle y_{i, v} + z_{v, i + 1} = 1, \forall v \in V, i = \overline{1, k - 1},
\end{array}
\label{order4}\\
& \begin{array}{l}
\displaystyle y_{i, u} + z_{u, i} +  y_{i, v} + z_{v, i} \ge 1, \forall uv \in E, i = \overline{1, k},
\end{array}
\label{order5}\\
& \begin{array}{l}
\displaystyle y_{i, q} - y_{i, v} \ge 0, \forall v \in V, i = \overline{1, k - 1},
\end{array}
\label{order6}\\
& \begin{array}{l}
\displaystyle y_{i, v}, z_{v, i} \in \{ 0, 1 \}, \forall v \in V, i = \overline{1, k},
\end{array} 
\label{orde71}
%\label{eqGlover8}
\end{align}
The vertex $ q $ will be assigned to the largest chosen color. 

\subsection{The Partial-Ordering Model\\ Adapted for Equitable Coloring}

In order to introduce the specific constraints for color classes cardinality we use the following result.
\begin{lemma}
\label{lemma1}
The color $ i \; ( \ge 2) $ class has cardinality
$ \displaystyle \left( \sum_{v \in V} y_{i-1, v} - \sum_{v \in V} y_{i, v} \right) $, while the color $ 1 $ class cardinality is $ \displaystyle \left( n - \sum_{v \in V} y_{1, v} \right) $. 
\end{lemma}
\proof{Each vector $ y_{\cdot, v} \in \{ 0, 1 \} ^k $ has its elements in increasing order:
\[ {y_{\cdot, v}}^T  = ( 1, 1, \ldots, 1, \underbracket{0}_{i}, \ldots, 0) \]
Vertex $ v $ receives the color $ i $ if and only if the first zero value occurs in position $ i $; this property allows us to compute the cardinalities of the color classes. The number of vertices having color at least $ i \; (\ge 2) $ is $ \displaystyle \sum_{v \in V} y_{i, v} $, while the number of vertices having color at least $ 1 $ is $ n $.}

Based on this lemma one can model the problem of deciding if a graph admits an equitable $ p $-coloring:
\begin{align}
(M1P) \:\: \: min & \begin{array}{l}
 \displaystyle 1
\end{array}
\notag\\
& \begin{array}{l}
\displaystyle z_{v, 1} =  0, \forall v \in V,
\end{array}
\label{eq_order1}\\
& \begin{array}{l}
\displaystyle y_{k, v} = 0, \forall v \in V,
\end{array}
\label{eq_order2}\\
& \begin{array}{l}
\displaystyle y_{i, v} - y_{i + 1, v} \ge  0, \forall v \in V, i = \overline{1, p - 1},
\end{array}
\label{eq_order3}\\
& \begin{array}{l}
\displaystyle y_{i, v} + z_{v, i + 1} = 1, \forall v \in V, i = \overline{1, p - 1},
\end{array}
\label{eq_order4}\\
& \begin{array}{l}
\displaystyle y_{i, u} + z_{u, i} +  y_{i, v} + z_{v, i} \ge 1, \forall uv \in E, i = \overline{1, p},
\end{array}
\label{eq_order5}\\
& \begin{array}{l}
\displaystyle \sum_{v \in V} y_{i - 1, v} - \sum_{v \in V} y_{i, v} - \ceil*{\frac{n}{p}} v_i -  \floor*{\frac{n}{p}} w_i = 0, \forall i = \overline{2, p},
\end{array}
\label{eq_order6}\\
& \begin{array}{l}
\displaystyle \sum_{v \in V} y_{1, v} + \ceil*{\frac{n}{p}} v_1 +  \floor*{\frac{n}{p}} w_1 = n,
\end{array}
\label{eq_order7}\\
& \begin{array}{l}
\displaystyle v_i + w_i = 1, \forall i = \overline{1, p},
\end{array}
\label{eq_order8}\\
& \begin{array}{l}
\displaystyle v_i, w_i, y_{i, v}, z_{v, i} \in \{ 0, 1 \}, \forall v \in V, i = \overline{1, p},
\end{array} 
\notag
%\label{eqGlover8}
\end{align}
Equations  \eqref{eq_order6} -- \eqref{eq_order8} constraining the cardinalities of the color classes can be replaced by (we get rid in this way of the variables $ v_i $ and $ w_i $, $ i = \overline{1, p} $):
\begin{align}
& \begin{array}{l}
\displaystyle \sum_{v \in V} y_{i - 1, v} - \sum_{j \in V} y_{i, v} \le \displaystyle \ceil*{\frac{n}{p}}, \forall i = \overline{2, p},
\end{array}
\tag{23'}\\
& \begin{array}{l}
\displaystyle \sum_{v \in V} y_{i - 1, v} - \sum_{v \in V} y_{i, v} \ge \displaystyle \floor*{\frac{n}{p}}, \forall i = \overline{2, p},
\end{array}
\tag{23''}\\
& \begin{array}{l}
\displaystyle n - \sum_{v \in V} y_{1, v} \le \displaystyle \ceil*{\frac{n}{p}},
\end{array}
\tag{24'}\\
& \begin{array}{l}
\displaystyle n - \sum_{v \in V} y_{1, v} \ge \displaystyle  \floor*{\frac{n}{p}},
\end{array}
\tag{24''}
\end{align}
One can observe that for these models there is no need of an objective function, the only question is if the subjacent polyhedra are non-empty. 

\subsection{The Maximum Cardinality of an Equitable Color Class}

The above models can be further modified for finding the maximum cardinality of a color class in an equitable coloring; such a parameter cannot give the equitable chromatic number but can give lower bounds for it. In the following we will suppose that $ k $ is an upper bound for $ \chi_{eq}(G) $ and $ M $ is a very large integer. The model is designed such that the maximum cardinality corresponds to the color $ 1 $.
\begin{align}
(M1) \:\: \:max & \begin{array}{l}
 \displaystyle \left( n- \sum_{v \in V} y_{1, v}\right) 
\end{array}
\notag\\
& \begin{array}{l}
\displaystyle z_{v, 1} =  0, \forall v \in V,
\end{array}
\label{eqN_order1}\\
& \begin{array}{l}
\displaystyle y_{k, v} = 0, \forall v \in V,
\end{array}
\label{eqN_order2}\\
& \begin{array}{l}
\displaystyle y_{i, v} - y_{i + 1, v} \ge  0, \forall v \in V, \forall i = \overline{1, k - 1},
\end{array}
\label{eqN_order3}\\
& \begin{array}{l}
\displaystyle y_{i, v} + z_{v, i + 1} = 1, \forall v \in V, \forall i = \overline{1, k - 1},
\end{array}
\label{eqN_order4}\\
& \begin{array}{l}
\displaystyle y_{i, u} + z_{u, i} +  y_{i, v} + z_{v, i} \ge 1, \forall uv \in E,  \forall i = \overline{1, k},
\end{array}
\label{eqN_order5}\\
& \begin{array}{l}
\displaystyle \sum_{v \in V} y_{i - 1, v} - \sum_{v \in V} y_{i, v} + \sum_{v \in V} y_{1, v} \le n, \forall i = \overline{2, k},
\end{array}
\label{eqN_order6}\\
& \begin{array}{l}
\displaystyle \sum_{v \in V} y_{i - 1, v} - \sum_{v \in V} y_{i, v} + \sum_{v \in V} y_{1, v} + M v_i \ge \displaystyle n -1,\forall i = \overline{2, k},
\end{array}
\label{eqN_order7}\\
& \begin{array}{l}
\displaystyle \sum_{v \in V} y_{i - 1, v} - \sum_{v \in V} y_{i, v} + M v_i \le M, \forall i = \overline{2, k},
\end{array}
\label{eqN_order8}\\
& \begin{array}{l}
\displaystyle v_i, y_{i, v}, z_{v, i} \in \{ 0, 1 \}, \forall v \in V, \forall i = \overline{1, k},
\end{array} 
\notag
%\label{eqGlover8}
\end{align}
Let $ u_i $ be the number of vertices having color $ i $, constraints \eqref{eqN_order6} prevent that no other color class has its cardinality strictly greater than $ u_1 $.  Constraints \eqref{eqN_order7} -- \eqref{eqN_order8} ensure that the cardinality of any color class is at least $ (u_1 - 1) $ or $ 0 $ (empty color classes cannot be avoided in this model), this can be done using the big $ M $ method:
\begin{align}
& \begin{array}{l}
 u_i \le M(1 - v_i)
\end{array}
\notag\\
& \begin{array}{l}
 u_i \ge u_1 - 1- M v_i
\end{array}
\notag
\end{align}
One may attempt to break the symmetry by requiring that the cardinalities are in non-increasing order: $ u_1 \ge u_2 \ge \ldots \ge u_k $, but experiments show no performance improvement by doing so. On the other hand the way in which we choose the very large integer $ M $ has a strong influence on solving the problem.

\begin{lemma}
\label{lemma2}
We could choose $ M  = \ceil*{n/k_0} $, where $ k_0 $ is a lower bound for $ \chi_{eq}(G) $.
\end{lemma}
\proof{ Suppose that the problem from above has a solution corresponding to a equitable $ p $-coloring of $ G $. If $ v_i $ is $ 0 $, then the first equation is satisfied if~$ M $ is an upper bound of $ u_i $ which is at most $ u_1 $ from \eqref{eqN_order6}, while if $ v_i $ is $ 1 $, then the second equation is satisfied if $ u_1 \le M + 1 $. Hence an appropriate value for $ M $ would be an upper bound for $ u_1 $ like $ \ceil*{n/k_0} $, since
\[ u_1 \le  \ceil*{n/p} \le \ceil*{n/\chi_{eq}(G)} \le \ceil*{n/k_0}.\] 

\vspace{-8ex}}

\subsection{The Assignment Model Revisited}

The assignment model can be modified for finding the maximum cardinality of an equitable color class (see \cite{mendez14b}).
\begin{align}
(M2) \: \: max & \begin{array}{l}
  \displaystyle \left( \sum_{v \in V} x_{v1} \right)
\end{array}
\notag\\
& \begin{array}{l}
\displaystyle \sum_{i = 1}^k x_{vi}=  1, \forall v \in V,
\end{array}
\label{asgnRev1}\\
& \begin{array}{l}
\displaystyle x_{ui} + x_{vi} - w_i \le 0, \forall uv \in E, \forall i = \overline{1, k}
\end{array}
\label{asgnRev2}\\
& \begin{array}{l}
\displaystyle w_{i +1} - w_i \le 0, \forall i = \overline{1, k - 1}
\end{array}
\label{asgnRev3}\\
& \begin{array}{l}
\displaystyle \sum_{v \in V} x_{vi} - \sum_{v \in V} x_{v1} \le 0, \; \forall i = \overline{2, k},
\end{array}
\label{asgnRev4}\\
& \begin{array}{l}
\displaystyle \sum_{v \in V} x_{vi} - \sum_{v \in V} x_{v1} + Mv_i \ge -1, \; \forall i = \overline{2, k},
\end{array}
\label{asgnRev5}\\
& \begin{array}{l}
\displaystyle \sum_{v \in V} x_{vi} + Mv_i \le M, \; \forall i = \overline{2, k},
\end{array}
\label{asgnRev6}\\
& \begin{array}{l}
\displaystyle x_{vi} \in \{ 0, 1 \} \;  \forall v \in V, \forall i = \overline{1, k}, \; v_i, w_i \in \{ 0, 1 \} \; \forall i = \overline{1, k}
\end{array}
\notag
\end{align}

\noindent Obviously, the same choice for $ M $ as above works here. For these reasons the implementations of both model $ M1 $ and $ M2 $ need a lower and an upper bound.

The corresponding model for deciding if the graph has an equitable $ k $-coloring follows
\begin{align}
(M2P) \: \: min  & 1
\notag\\
& \begin{array}{l}
\displaystyle \sum_{i = 1}^p x_{vi}=  1, \forall v \in V,
\end{array}
\label{asgnRevP1}\\
& \begin{array}{l}
\displaystyle x_{ui} + x_{vi} - w_i \le 0, \forall uv \in E, \forall i = \overline{1, p}
\end{array}
\label{asgnRevP2}\\
& \begin{array}{l}
\displaystyle w_{i +1} - w_i \le 0, \forall i = \overline{1, p - 1}
\end{array}
\label{asgnRevP3}\\
& \begin{array}{l}
\displaystyle \sum_{v \in V} x_{vi}  \le  \ceil*{\frac{n}{p}}, \; \forall i = \overline{1, p},
\end{array}
\label{asgnRevP4}\\
& \begin{array}{l}
\displaystyle \sum_{v \in V} x_{vi}  \ge  \floor*{\frac{n}{p}}, \; \forall i = \overline{1, p},
\end{array}
\label{asgnRevP5}\\
& \begin{array}{l}
\displaystyle x_{vi} \in \{ 0, 1 \} \;  \forall v \in V, \forall i = \overline{1, p}, \; w_i \in \{ 0, 1 \} \; \forall i = \overline{1, p}
\end{array}
\notag
\end{align}

\subsection{Lower Bound for Equitable Chromatic Number}

The following result shows that one may improve the lower bound for the equitable chromatic number by just providing upper bounds for the optimum in the above ILP problems ($ M1 $ and $ M2 $).

\begin{lemma}
Let $ G = (V, E) $ be a graph and $ \beta_0 $ be an integer upper bound for the maximum cardinality of a color class in any equitable colorings of the vertices of $ G $. Then $ \chi_{eq}(G) \ge \displaystyle \ceil*{\frac{n}{\beta_0}} $.

\end{lemma}
\label{lemma3}
\proof{Let $ k $ be the number of colors for an equitable coloring for which the maximum cardinality of a color class is $ \beta $, obviously
\[ \ceil*{\frac{n}{\beta}} \le k \le \ceil*{\frac{n}{\beta - 1}}. \]
Consider now an optimum equitable coloring and the maximum cardinality of one of its a color classes $ \beta_{eq} $. If $ \beta_{eq} < \beta $, then 
\[ \ceil*{\frac{n}{\beta_0}} \le \ceil*{\frac{n}{\beta}} \le \ceil*{\frac{n}{\beta_{eq}}} \le \chi_{eq}(G) \le \ceil*{\frac{n}{\beta_{eq} - 1}}. \]
} 

The proof of the following consequence is obvious.
\begin{corollary}
Suppose that $ \beta_1 $ is the maximum cardinality of a color class (i. e., a solution to the above ILP problem), then the maximum cardinality of a color class in an optimum maximum cardinality of a color class, $  \beta_{eq} $, belongs to the set
\[ \left \{ \beta \in {\mathbb Z}^* \: : \: \ceil*{\frac{n}{\beta - 1}} = \ceil*{\frac{n}{\beta_1}} \right\}. \]
\end{corollary}

\section{Numerical Experiments and Conclusions}\label{section:numerical}

In this section we present and analyze computational experiments all carried on an Intel i5-7500 CPU\@ 3.40 GHz with with 8 GB of memory on Ubuntu 18.04.5 LTS and using a Gurobi Academic License (Benchmarks~\cite{benchmark02} user time: r500.5=4.74 s). The benchmark instances (commonly used in the literature for the classic graph coloring problem) are available at: 
\url{http://cedric.cnam.fr/~porumbed/graphs/} and\\ \url{https://mat.tepper.cmu.edu/COLOR02/INSTANCES/}

\subsection{Method}
Our models are used in the following way for improving (if possible) lower bounds for the equitable chromatic number:

\begin{itemize}

\item[1.] Use first one of the two models $ (M1) $ or $ (M2) $ for finding upper bounds of the maximum cardinality of an equitable color class. 

\item[2.] If the found lower bounds are better update them and then iteratively try to improve them by employing model $ (M2) $ or $ (M2P) $.

\end{itemize}

\noindent In the first step from above we need lower and upper bounds for the equitable chromatic number; one can use $ \Delta(G) + 1 $ as an upper bound (\cite{hajnal70,kierstead10}) and~$ 3 $ as a lower bound - since the graphs are non-bipartite. In our experiments we used the best known bounds from the literature (\cite{mendez14a,lai15}), except for the  queen$ k $\_$ k $ graphs for which $ \Delta(G) + 1 $ was used as upper bound. It is interesting to note that it's not necessarily to completely solve models $ (M1) $ and ($ M2 $); using the Gurobi solver one can get an upper bound for the objective and take its integer part as the desired upper bound of the maximum cardinality of an equitable color class. 

In the second step, starting with  the (new discovered) lower bound we verify if the graph has an equitable coloring with the corresponding number of colors and increase this bound if the problem proves to be infeasible. This time we use the models  ($ M2 $) and ($ M2P $); ($ M2 $) can be used in the following way to verify if the given graph has an equitable $ p $-coloring: take $ M = \ceil*{n/p} $ - see Lemma \ref{lemma3} - and $ k = p $. Since the corresponding results for ($ M1 $) and ($ M1P $) are much less encouraging  we skipped them; our results show that the partial ordering base models are weaker than the assignment based models at least in terms of the continuous LP-relaxation.

\begin{table}[!h]
\caption{Numerical results for step 1.}
\footnotesize
\begin{tabular}{llllllllllll}%{p{0.11\textwidth}p{0.037\textwidth}p{0.037\textwidth}p{0.037\textwidth}p{0.037\textwidth}p{0.037\textwidth}}
\hline
\multirow{2}{*}{Instance} & \multirow{2}{*}{$ n $} & \multirow{2}{*}{$ m $} & \multirow{2}{*}{LB} & \multirow{2}{*}{UB} & \multicolumn{3}{c}{$ M1 $} & & \multicolumn{3}{c}{$ M2 $}   \\
\cline{6-8} \cline{10-12}
  & & & & & $ \beta_0 $ & LB* & time(s) &  & $ \beta_0 $ & LB* & time(s) \\
\hline
R125.1 & $ 125 $ & $ 209 $ & $ 3 $ & $ 5 $ & $ 25 $ & $ \mathbf{5} $ & $ 2.91^i $ &  & $ 30 $ & $ \mathbf{5} $ &  $ 0.04^i $  \\
R125.5 & $ 125 $ & $ 3838 $ & $ 3 $ & $ 36 $ & $ 5 $ & $ \mathbf{25} $ & $ 2.52^r $ & &  $ 5 $ & $ \mathbf{25} $  & $ 4.16^r $  \\
R250.1 & $ 250 $ & $ 867 $ & $ 3 $ & $ 8 $ & $ 70 $ & $ \mathbf{4} $ & $ 0.06^r $ & &  $ 36 $ & $ \mathbf{7} $ & $ 0.47^i $  \\
R250.5 & $ 250 $ & $ 14849 $ & $ 3 $ & $ 66 $ & $ 7 $ & $ \mathbf{36} $ & $ 261^b $&  & $ 7 $ & $ \mathbf{36} $ & $ 512^b$  \\
le450\_5c & $ 450 $ & $ 9803 $ & $ 3 $ & $ 5 $ & $ 93 $ & $ \mathbf{5} $ & $ 8.94^r $ &  & $ 92 $ & $ \mathbf{5} $ & $ 388.85^r $  \\
le450\_15c & $ 450 $ & $ 16680 $ & $ 3 $ & $ 15 $ & $ 49 $ & $ \mathbf{10}  $ & $ 22.78^r $ &  & $ 31 $ & $ \mathbf{15} $ & $ 222.77^r $  \\
le450\_25c & $ 450 $ & $ 17343 $ & $ 3 $ & $ 26 $ & $ 49 $ & $ \mathbf{10} $ & $ 91^b $ & &  $ 21 $ & $ \mathbf{22} $ & $ 51.86^r $  \\
le450\_25d & $ 450 $ & $ 17425 $ & $ 25 $ & $ 26 $ & $ 19 $ & $ 24 $ & $ 11.89^r $ & &  $ 18  $ & $ 25 $ & $ 31.29^r $ \\
flat300\_28\_0 & $ 300 $ & $ 21695 $ & $ 11 $ & $ 34 $ & $ 29 $ & $ 11 $ & $ 71.35^r $ &  & $ 28 $ & $ 11 $ & $ 194.05^r $  \\
flat300\_20\_0 & $ 300 $ & $ 21375 $ & $ 11 $ & $ 34 $ & $ 29 $ & $ 11 $ & $ 69.68^r $ & &  $ 28 $ & $ 11 $ & $ 140.71^r $ \\
ash608GPIA & $ 1216 $ & $ 7844 $ & $ 3 $ & $ 4 $ & $ 407 $ & $ 3 $ & $ 15.31^r $ & &  $ 406 $ & $ 3 $ & $ 30.71^r $  \\
ash958GPIA & $ 1916 $ & $ 12506 $ & $ 3 $ & $ 4 $ & $ 640 $ & $ 3 $ & $ 6.88^r $&  & $ 639 $ & $ 3 $ & $ 84.54^r $ \\
DSJC125.5 & $ 125 $ & $ 7782 $ & $ 9 $ & $ 17 $ & $ 15 $ & $ 9 $ & $ 1.12^r $ &  & $ 14 $ & $ 9 $ & $ 7.15^r $ \\
DSJC125.9 & $ 125 $ & $ 13922 $ & $ 43 $ & $ 44 $ & $ 4 $ & $ 32 $ & $ 10.43^r $ &  & $ 3 $ & $ 42 $ & $ 1.82^r $ \\
DSJC250.1 & $ 250 $ & $ 6436 $ & $ 4 $ & $ 8 $ & $ 64 $ & $ 4 $ & $1.32^r $&  & $ 63 $ & $ 4 $ & $ 15.00^r $ \\
DSJC250.5 & $ 250 $ & $ 31336 $ & $ 12 $ & $ 30 $ & $ 21 $ & $ 12 $ & $34.75^r $ & & $ 21 $ & $ 12 $ & $ 23.63^r $  \\
inithx.i.3 & $ 621 $ & $ 13969 $ & $ 3 $ & $ 37 $ & $ 208 $ & $ 3 $ & $ 8.60^r $&  & $ 22 $ & $ \mathbf{29} $ & $ 283^b $   \\
inithx.i.2 & $ 645 $ & $ 13979 $ & $ 30 $ & $ 36 $ & $ 23 $ & $ 29 $ & $ 40.16^r $ &  & $ 22 $ & $ 30 $ & $ 20.50^r $ \\
mulsol.i.2 & $ 188 $ & $ 3885 $ & $ 34 $ & $ 36 $ & $ 7 $ & $ 27 $ & $ 5.55^i $ & & $ 6 $ & $ 32 $ & $ 10.56^i $ \\
2-Insertions\_5 & $ 597 $ & $ 3936 $ & $ 3 $ & $ 6 $ & $ 200 $ & $ 3 $ & $ 0.39^r $ &  & $ 199 $ & $ 3 $ & $ 3.78^r $  \\
1-Insertions\_6 & $ 607 $ & $ 6337 $ & $ 3 $ & $ 7 $ & $ 204 $ & $ 3 $ & $ 1.03^r $ &  & $ 203 $ & $ 3 $ & $ 4.27^r $  \\
4-FullIns\_4 & $ 690 $ & $ 6650 $ & $ 6 $ & $ 8 $ & $ 116 $ & $ 6 $ & $ 5.61^r $ &  & $ 115 $ & $ 6 $ & $ 4.20^r $  \\
4-FullIns\_5 & $ 4146 $ & $ 77305 $ & $ 6 $ & $ 9 $ & $ 692 $ & $ 6 $ & $ 55.65^r $ & & $ 691 $ & $ 6 $ & $ 604.78^r $  \\
wap02a & $ 2464 $ & $ 111742 $ & $ 40 $ & $ 41 $ & $ - $ & $ - $ & $ - $ & & $ 62 $ & $ 40 $ & $ 549.10^r $   \\
wap05a & $ 905 $ & $ 43081 $ & $ 3 $ & $ 50 $ & $ - $ & $ - $ & $ - $  & & $ 46 $ & $ \mathbf{19} $ & $ 211.53^r $  \\
wap06a & $ 947 $ & $ 43571 $ & $ 3 $ & $ 41 $ & $ 85 $ & $ \mathbf{12} $ & $ 62.84^r $ & & $ 31 $ & $ \mathbf{31}$ & $ 160.85^r $\\
wap07a & $ 1809 $ & $ 103368 $ & $ 3 $ & $ 43 $ & $ - $ & $ - $ & $ - $ & & $ 58 $ &  $ \mathbf{32} $ & $ 1667.27^r $ \\
wap08a & $ 1870 $ & $ 104176 $ & $ 3 $ & $ 43 $ & $ - $ & $ - $ & $ - $ & & $ 68 $ & $ \mathbf{28} $ & $ 903.94^r $  \\
queen12\_12 & $ 144 $ & $ 5192 $ & $ 3 $ & $ 44 $ & $ 12 $ & $ \mathbf{12} $ & $ 2.23^r $ & & 12 & $ \mathbf{12} $ &  $ 6.06^r $ \\
queen13\_13 & $ 144 $ & $ 6656 $ & $ 13 $ & $ 48 $ & $ 13 $ & $ \mathbf{13} $ & $ 5.16^r $ & & $ 13 $ & $ \mathbf{13} $ & $ 59.23^r $  \\
queen14\_14 & $ 196 $ & $ 8372 $ & $ 3 $ & $ 51 $ & $ 14 $ & $ \mathbf{14} $ & $ 2.33^r $ & & $ 14 $ & $ \mathbf{14} $ & $ 7.36^r $  \\
queen15\_15 & $ 225 $ & $ 10360$ & $ 3 $ & $ 56 $ & $ 15 $ & $ \mathbf{15} $ & $ 4.75^r $ & & $ 15 $ & $ \mathbf{15} $ & $ 16.19^r $  \\
queen16\_16 & $ 256 $ & $ 12640$ & $ 3 $ & $ 59 $ & $ 16 $ & $ \mathbf{16} $ & $ 7.40^r $ & & $ 17 $ & $ \mathbf{16} $ &  $ 24.03^r $ \\
\hline
\end{tabular}\label{table1}
\end{table}

The computational results of the first step are reported in Table \ref{table1}. Columns 1, 2, and 3 display the name of the instance, the number of vertices, and the number of edges, correspondingly; columns 4 and 5 show the known bounds for the equitable chromatic number. Columns 6-8 and 9-11 display the upper bound for the maximum cardinality of an equitable color class (called $ \beta_0 $), the corresponding lower bound for $ \chi_{eq} $, and the time needed to find this upper bound, or "$ - $" if the solver was not able to find it within the~$ 30 $ minutes time limit.

There are three ways of finding $ \beta_0 $ using Gurobi solver: by finding an optimal solution to the corresponding ILP (marked by an "$ i $"), by computing the root relaxation objective ("$ r $"), or by finding upper bounds of the ILP objective function ("$ b $").

For the thirty-three evaluated instances we found seventeen new lower bounds, in only two cases our results were not better, while for the remaining instances the best known bounds were confirmed. The $ M2 $ model proved to give better bounds but needs more time for almost all instances, while the~$ M1 $ model is faster but failed to deliver any bound in four cases.

\begin{table}[!h]
\caption{Numerical results for steps 1 and 2.}
\footnotesize
\begin{tabular}{llllllll}
\hline
\multirow{2}{*}{Instance} & new (old)  & \multicolumn{2}{c}{$ M2 $} & & \multicolumn{2}{c}{$ M2P $}  & \multirow{2}{*}{Remarks}\\
\cline{3-4} \cline{6-7}
  & LB &  LB* & time(s) & &  LB* & time(s) &   \\
\hline
R125.1 & $ \mathbf{5} $ ($ 3 $)  & $ - $ & $ - $ & & $ - $ & $ - $ & step 1, $ \chi_{eq} = \mathbf{5} $ \\
R125.5 & $ \mathbf{35} $ ($ 3 $) & $ 35 $ & $ 18.10 $ & & $ 35 $ & $ 2.01 $ &  step 2 \\
R250.1 & $ \mathbf{8} $ ($ 3 $) & $ 8 $ & $ 0.07 $ & & $ 8 $ & $ 0.03 $ & step 2, $ \chi_{eq} = \mathbf{8} $  \\
R250.5 & $ \mathbf{57} $ ($ 3 $) & $ 57 $ & $ 36.62 $ & & $ 57 $ & $ 31.39 $ & step 2 \\
le450\_5c & $ \mathbf{5} $ ($ 3 $) & $ - $ & $ - $ & & $ - $ & $ - $ & step 1, $ \chi_{eq} = \mathbf{5} $ \\
le450\_15c & $ \mathbf{15} $ ($ 3 $)  & $ - $ & $ - $ & & $ - $ & $ - $ & step 1, $ \chi_{eq} = \mathbf{15} $ \\
le450\_25c & $ \mathbf{25} $ ($ 3 $) & $ 25 $ & $ 61.22 $ & & $ 25 $ & $ 16.60 $ & step 2  \\
le450\_25d & $ 25 $ ($ 25 $) & $ - $ & $ - $ & & $ - $ & $ - $ & step 1  \\
flat300\_28\_0 & $ 11 $ ($ 11 $) & $ - $ & $ - $ & & $ - $ & $ - $ & step 1  \\
flat300\_20\_0 & $ 11 $ ($ 11 $) & $ - $ & $ - $ & & $ - $ & $ - $ & step 1  \\
ash608GPIA & $ \mathbf{4} $ ($ 3 $) & $ 4 $ & $ 120.10 $ & & $ 4 $ & $ 0.13 $ & step 2, $ \chi_{eq} = \mathbf{4} $ \\
ash958GPIA & $ \mathbf{4} $ ($ 3 $) & $ 4 $ & $ 322.02 $ & & $ 4 $ & $ 0.22 $ & step 2, $ \chi_{eq} = \mathbf{4} $ \\
DSJC125.5 & $ 9 $ ($ 9 $) & $ - $ & $ - $ & & $ - $ & $ - $ & step 1  \\
DSJC125.9 & $ 42 $ ($ 43 $) & $ - $ & $ - $ & & $ - $ & $ - $ & step 1  \\
DSJC250.1 & $ \mathbf{5} $ ($ 4 $) & $ 5 $ & $ 20.10 $ & & $ 5 $ & $ 7.18 $ & step 2 \\
DSJC250.5 & $ 12 $ ($ 12 $) & $ - $ & $ - $ & & $ - $ & $ - $ &  step 1 \\
inithx.i.3 & $ \mathbf{31} $ ($ 3 $) & $ 31 $ & $ 32.34 $ & & $ 31 $ & $ 21.08 $ & step 2 \\
inithx.i.2 & $ \mathbf{31}  $ ($ 30 $) & $ 31 $ & $ 28.51 $ & & $ 31 $ & $ 14.38 $ & step 2 \\
mulsol.i.2 & $ 32 $ ($ 34 $) & $ - $ & $ - $ & & $ - $ & $ - $ & step 1 \\
2-Insertions\_5 &  $ \mathbf{4} $ ($ 3 $) & $ 4 $ & $ 91.20 $ & & $ 4 $ & $ 56.23 $ & step 2 \\
1-Insertions\_6 & $ \mathbf{4} $ ($ 3 $) & $ 4 $ & $ 42.74 $ & & $ 4 $ & $ 6.05 $ & step 2 \\
4-FullIns\_4 & $ \mathbf{7} $ ($ 6 $) & $ 7 $ & $ 36.69 $ & & $ 7 $ & $ 17.90 $ & step 2 \\
4-FullIns\_5 & $ 6 $  ($ 6 $) & $ - $ & $ - $ & & $ - $ & $ - $ & step 1\\
wap02a & $ 40 $ ($ 40 $) & $ - $ & $ - $ & & $ - $ & $ - $ & step 1 \\
wap05a & $ \mathbf{46} $ ($ 3 $) & $ 46 $ & $ 33.84 $ & & $ 46 $ & $ 26.45 $ & step 2 \\
wap06a & $ \mathbf{40} $ ($ 3 $) & $ 40 $ & $ 47,64 $ & & $ 40 $ & $ 33.33 $ & step 2 \\
wap07a & $ \mathbf{40} $ ($ 3 $) & $ 40 $ & $ 442,90 $ & & $ 40 $ & $ 341.74 $ & step 2  \\
wap08a & $ \mathbf{40} $  ($ 3 $)  & $ 40 $ & $ 311.38 $ & & $ 40 $ & $ 283.97 $ & step 2 \\
queenk\_k & $ \mathbf{k} $ ($ 3 $) & $ - $ & $ - $ & & $ - $ & $ - $ & step 1 \\
\hline
\end{tabular}\label{table2}
\end{table}

Table \ref{table2} reports the numerical results achieved in the second step and outlines both steps. Column 1 display the name of the instance, column 2 shows the new lower bound found after performing both steps (in parentheses we have the older best known bounds) - if $ \chi_{eq} $ was already found after the first step, or if the second step failed to improve the first step, then here we have the result from Table \ref{table1}. Column 3-4 and 5-6 report the lower bound obtained after performing the second step (if any) and the average time. 

This step is a sequential procedure (for this reason we reported the average time), e. g. for R250.5.col instance the known lower bound and upper bound are $ 3 $ and $ 66 $, respectively, the first step enlarges the lower bound to $ 36 $, in the second step we verify if the graph admits equitable colorings with $ k = 36, 37, \ldots $ - and it turned out that the corresponding ILP problems are infeasible up to $ k = 56 $, hence $ \chi_{eq} \ge 57 $. When the solver finds integer solutions to the ILP, or the lower bound  $ \chi_{eq} $ is found and reported.

Column 7 contains the remarks concerning the decisive (and final) step and the finding of $ \chi_{eq} $. For twenty-four of the instances new lower bounds were found (an increase from the first step) and for six of them $ \chi_{eq} $ was found (three are due only to the first step).

Model $ M2P$ is faster than $ M2 $, but both models gave same bounds - or both failed to improve them; the same $ 30 $ minutes time limit - for each iteration - was used in this step also.

\section{Conclusions}\label{section:conclusion}

In this paper we introduce and analyze a two stage procedure for improving the lower bound for the equitable chromatic number of a graph. This is a new approach since usually the literature knows only procedures that decreases the upper bound for this number.

Our approach is based on finding the maximum cardinality (or, at least upper bound) of an equitable color class. Our method finds improved lower bounds  for $ 25 $ out of $ 33 $ investigated instances.

The first step of our method employs a new model - based on a partial ordering model for coloring - which proved to be effective but less efficient than the classic assignment model, while the second step is based only on the classic assignment model. Both models could be subject to improvements using cuts and this will be a new line of research especially for the (newer) partial ordering model.

\bibliography{tr-template_olariu}
% if you use BibTeX then
% \bibliography{BibTeX file}
% else

%\newpage
%
%\hypertarget{lastpage}{}
%\lfoot{{\bf \copyright\hspace{0.01mm} Scientific Annals of Computer Science yyyy}}
\end{document}